\documentclass[12pt]{article}

\usepackage{amsmath}
\usepackage{amsthm}
\usepackage{amsfonts}
\usepackage{amssymb}
\usepackage[all]{xy}
\usepackage{amsmath,amscd}

\theoremstyle{plain}
\newtheorem{theorem}{Theorem}[section]
\newtheorem{proposition}{Proposition}[section]
\newtheorem{lemma}{Lemma}[section]

\theoremstyle{definition}

\newcommand{\Q}{\mathbb{Q}}
\newcommand{\R}{\mathbb{R}}
\newcommand{\Z}{\mathbb{Z}}

\newcommand{\PP}{\mathbb{P}}
\newcommand{\RR}{\mathbb{R^{+}}}

\newcommand{\Bl}{\operatorname{Bl}}

\newcommand{\NNE}{\overline{\operatorname{NE}}}
\newcommand{\Nef}{\operatorname{Nef}}

\newcommand{\Eff}{\operatorname{Eff}}

\newcommand{\mult}{\operatorname{mult}}

\newcommand{\Exc}{\operatorname{Exc}}

\newcommand{\Supp}{\operatorname{Supp}}
\newcommand{\OO}{{\cal O}}

\newcommand{\wt}{\widetilde}

\newcommand{\st}{_{*}^{-1}}

\title
{Some examples of log Fano structures on blow-ups along subvarieties in products of two projective spaces }
\author{Toru Tsukioka}

\begin{document}

\maketitle 

\begin{abstract}
We consider the problem to determine which blow-ups along subvarieties in products of two projective spaces are log Fano. By describing the nef cones of such blow-ups with special centers, we give a partial classification result. For each example, explicit boundary divisors are also given.  
\end{abstract}



\section{Introduction}
In the study of higher dimensional Fano varieties, it is interesting to find examples 
with special geometric structures. 
The recent papers \cite{AM} and \cite{LP} investigate log Fano manifolds obtained by blowing up points in projective spaces  
or products of projective spaces. 
In \cite{Tsukioka}, the author considered weak Fano manifolds 
having small contractions (i.e. birational morphisms whose exceptional loci has codimension greater than or equal to 2). 

We recall the construction of a small contraction using successive blow-ups. Let $Y$ be a smooth projective variety of dimension $n \geq 3$. Let $S$ and $C$ be smooth subvarieties of codimension $\geq 2$ 
with $\dim S +\dim C\leq n-1$. Assume that $S$ and $C$ intersect transversally at points. Let $X$ be the blow-up along $C$ and $\wt{X}$ the blow-up along the strict transform of $S$. Then, $\wt{X}$ has exceptional loci of codimension $\geq 2$ over the intersections of $S$ and $C$ (this is a slight generalization of \cite{Kawamata} Example 2.6). 

The problem is to determine the triples $(Y,S,C)$ such that $\wt{X}$ is (weak) Fano or more generally log Fano. In this note, we consider the case in which $Y$ is a product of two projective spaces, 
$S$ is a fiber of a projection, and $C$ is a hypersurface in a fiber of the other projection.

We work over the field of complex numbers.

\begin{theorem} Let $Y=\PP^{n-k}\times \PP^{k}$ with $n\geq 3$ and $k\in \{2,\cdots, n-1\}$.
Let $S$ be a fiber of the projection $Y\to \PP^{k}$. Let $C$ be a smooth hypersurface of degree 
$d\geq 1$ in a fiber of the other projection $Y\to \PP^{n-k}$. Assume $S\cap C\neq \emptyset$.  
Let $\pi:X\to Y$ be the blow-up along $C$. Let $S'$ be the strict transform of $S$ by $\pi$. 
Let $\wt{X}$ be the blow-up of $X$ along $S'$. Then, the anti-canonical divisor of $\wt{X}$ is big. 

If $n=3$, $\wt{X}$ is log Fano for any $d\geq 1$, 
is weak Fano for $d\in\{1,2,3\}$, and is not Fano for any $d\geq 1$. 

Assume $n\geq 4$. 
If one of the following holds, then $\wt{X}$ is log Fano. 
\begin{itemize}
\item[(1)] $d=1$, $n$ is even and $k=\frac{n}{2}$
\item[(2)] $d=1$, $n$ is even and $k=\frac{n+2}{2}$
\item[(3)] $d=1$, $n$ is odd and $k\in\{\frac{n-1}{2},\frac{n+1}{2}\}$
\item[(4)] $(n,k,d)=(4,2,2)$ or $(4,3,2)$
\item[(5)] $(n,k,d)=(5,3,2)$. 
\end{itemize}
Moreover, $\wt{X}$ is Fano if and only if (1) holds, and 
$\wt{X}$ is weak Fano if and only if (1),(3) or (5) holds.
\end{theorem}
Remark: In Section \ref{proof}, we give explicit log Fano pairs for the cases (1)--(5). We also find a condition for divisors in a specific form to be boundary divisors (Lemma \ref{log Fano pair}). The author do not know whether $\wt{X}$ is log Fano or not in the other cases.

\

Definition: Let $X$ be a smooth projective variety. 
We call $X$ {\it Fano} (resp. {\it weak Fano}) if the anti-canonical divisor $-K_{X}$ is ample (resp. nef and big). 
A {\it log Fano pair} $(X,\Delta)$ is a klt pair such that $-(K_{X}+\Delta)$ is ample. 
We call $X$ {\it log Fano} if there exists an effective $\Q$-divisor $\Delta$ such that $(X,\Delta)$ is a log Fano pair. The divisor $\Delta$ is called {\it boundary divisor}. We have the following implications: 
\begin{center}
Fano $\Rightarrow$ weak Fano $\Rightarrow$ log Fano $\Rightarrow$ $-K_{X}$ is big.
\end{center}
The first implication is obvious. For the second one, see \cite{AM} Lemma 2.5. 
The last one is due to Kodaira's lemma (see \cite{Lazarsfeldbook} Corollary 2.2.7). 

\

This note is organized as follows. 
First, we find a nef divisor on a special prime divisor $\wt{D}\subset \wt{X}$ 
(Section \ref{divisor}).
Then, we determine the structure of the nef cone of $\wt{X}$ 
by using the restriction map $\Nef (\wt{X})\to \Nef (\wt{D})$. 
We also compute the effective cone of $\wt{X}$ (Section \ref{cone}).  
Once the nef cone of $\wt{X}$ is described explicitly, it is easy to find conditions for $-(K_{\wt{X}}+\Delta)$ to be ample. 
We prove Theorem by solving a system of linear inequalities in the coefficients of the divisor $\Delta$ and 
the values $n$, $k$ and $d$ (Section \ref{proof}). 

\

The numerical equivalence class of an $\R$-divisor $D$ is denoted by $[D]$. 
Let $\RR$ be the set of non-negative real numbers. We put 
$$\RR[D_{1},\cdots,D_{k}]:=\RR[D_{1}]+\cdots +\RR[D_{k}]. $$

\section{A nef divisor on a special prime divisor}\label{divisor}

The following is useful to determine the structure of 
nef cones.

\begin{lemma} (\cite{Tsukioka} Lemma 3.2)\label{induction}
Let $X$ be a smooth projective variety and 
let $D$ be a smooth prime divisor on $X$. 
Let $N$ be a divisor on $X$. 
If $[N] |_{D} \in \Nef (D)$ and $[N-D]\in \Nef (\wt{X})$,  
then $[N]\in \Nef (X)$.  
\end{lemma}
\begin{lemma}\label{point}
Let $\alpha: X\to \PP^{n}$ ($n\geq 2$) be the blow-up at one point. 
Then, we have 
\begin{eqnarray*}
\Nef (X)=\RR[\alpha^{*}\OO_{\PP^{n}}(1),\alpha^{*}\OO_{\PP^{n}}(1)-\Exc (\alpha)]. 
\end{eqnarray*} 
In particular, $\alpha^{*}\OO_{\PP^{n}}(1)-\Exc (\alpha)$ is nef. 
\end{lemma}
{\it Proof.} Note that we have $\rho(X)=2$. Hence, the cone $\Nef (X)$ has exactly two extremal rays. Hence, the statement follows from 
the fact that the linear system 
$|\alpha^{*}\OO_{\PP^{n}}(1)-\Exc (\alpha)|$ is base point free. \qed

\

We fix the following notation (valid only in this section): 

Let $m$ and $k$ be natural numbers such that $m\geq 1$ and $k\geq 2$. 
Let $V$ be a smooth hypersurface of degree $d\geq 1$ in $\PP^{k}$. 
Put $D:=\PP^{m}\times V$.  
Let $pr_{1}: D\to \PP^{m}$ and $pr_{2}: D\to V$ be the projections. 
Let $C$ (resp. $S$) be a fiber of $pr_{1}$ (resp. $pr_{2}$). 
Note that $C$ (resp. $S$) is isomorphic to $V$ (resp. $\PP^{m}$). 
If $m=1$ (resp. $k=2$), then $C$ (resp. $S$) is a divisor on $D$.  
Put $H:=pr_{1}^{*}\OO_{\PP^{m}}(1)$
and 
$L:=pr_{2}^{*}\OO_{V}(1)$ where $\OO_{V}(1):=\OO_{\PP^{k}}(1)|_{V}$. 
We define the morphisms $\pi:D'\to D$ and $\beta:\wt{D}\to D'$ as follows.
\begin{itemize}
\item If $C$ in not a divisor on $D$, 
let $\pi:D'\to D$ be the blow-up along $C$. 
We put 
$E:=\Exc (\pi)$, $H':=\pi^{*}H$, $L':=\pi^{*}L$ and $S':=\pi\st S$.
If $C$ is a divisor on $D$, 
we put
$\pi:=id$, $D':=D$, $E:=C$, $H':=H$, $L':=L$ and  $S'=S$.  
\item If $S'$ is not a divisor on $D'$, 
let $\beta :\wt{D}\to D'$ be the blow-up along $S'$. We put 
$F:=\Exc (\beta)$, $\wt{H}:=\beta^{*}H'$, $\wt{L}:=\beta^{*}L'$ and  $\wt{E}:=\beta^{*}E$. 
If $S'$ is a divisor on $D'$, we put 
$\beta:=id$, $\wt{D}:=D'$, $F:=S'$, $\wt{H}:=H'$, $\wt{L}:=L'$ and  $\wt{E}:=E$. 
\end{itemize}

\begin{proposition}\label{nef} We have 
$[\wt{H}+d\wt{L}-\wt{E}-F]\in \Nef(\wt{D}).$
\end{proposition}
{\it Proof}. Since  $\wt{H}+d\wt{L}-\wt{E}-F=(d-1)\wt{L}+(\wt{H}-\wt{E})+(\wt{L}-\wt{F})$, 
it is sufficient to show that  $\wt{H}-\wt{E}$ and $\wt{L}-F$ are nef. 

Case $m=1$: 
$C$ is a divisor on $D=\PP^{1}\times V$ and 
we have $H'=H=pr_{1}^{*}\OO_{\PP^{1}}(1)\sim C=E$. 
Hence, $\wt{H}-\wt{E}=\beta^{*}(H'-E)$ is linearly equivalent to $0$, in particular nef.
Consider the blow-up 
$\mu:\Bl_{s}(\PP^{k})\to \PP^{k}$ at the point 
$s:=pr_{2}(S)$. Put $\wt{V}:=\mu\st V$. 
Then, $\mu |_{\wt{V}}:\wt{V}\to V$ is the blow-up at the point $s\in V$. 
Let $\wt{pr_{2}}:\wt{D}=\PP^{1}\times \wt{V}\to\wt{V}$ be the projection.  
Then, we have the commutative diagram: 
\begin{eqnarray*}
\xymatrix{
\wt{D}\ar[d]_{\wt{pr_{2}}} \ar[r]^{\beta} & D \ar[d]^{pr_{2}} \\ 
\wt{V}\ar[r]_{\mu |_{\wt{V}}}              & V.
}
\end{eqnarray*}
This yields 
\begin{eqnarray*}
\wt{L}&=&\beta^{*}L
=\beta^{*}pr_{2}^{*}\OO_{V}(1)
=\wt{pr_{2}}^{*}\left(\mu^{*}\OO_{\PP^{k}}(1)|_{\wt{V}}\right), \\ 
F&=&\wt{pr_{2}}^{*}\Exc (\mu |_{\wt{V}})
=\wt{pr_{2}}^{*}\left(\Exc (\mu)|_{\wt{V}}\right). 
\end{eqnarray*}
By Lemma \ref{point}, $\mu^{*}\OO_{\PP^{k}}(1)-\Exc (\mu)$ is nef. 
Hence, so is $\wt{L}-F$. 

Case $m\geq 2$: 
Consider the blow-up $\varepsilon: \Bl_{c}(\PP^{m})\to \PP^{m}$ at the point $c:=pr_{1}(C)$. 
Let $pr_{1}':D'\simeq \Bl_{c}(\PP^{m})\times V\to \Bl_{c}(\PP^{m})$ be the projection. 
Then, we have the commutative diagram: 
\begin{eqnarray*}
\xymatrix{
D' \ar[d]_{\pi} \ar[r]^{pr_{1}'} & \Bl_{c}(\PP^{m}) \ar[d]^{\varepsilon} \\ 
D \ar[r]_{pr_{1}} & \PP^{m}.}
\end{eqnarray*}
By Lemma \ref{point} again, 
$H'-E=(pr_{1}')^{*}(\varepsilon^{*}\OO_{\PP^{m}}(1)-\Exc (\varepsilon))$ 
is nef. Thus, so is $\wt{H}-\wt{E}=\beta^{*}(H'-E)$. 

Now, we show $\wt{L}-F$ is nef. Fix $m\geq 2$ and use the induction on $k\geq 2$. 
If $k=2$, then $\dim V=1$ and $S$ is a divisor on $D$. 
Since $S$ is a fiber of $pr_{2}:D\to V$, we have 
$L=pr_{2}^{*}(\OO_{\PP^{k}}(1) |_{V})\equiv dS$ and $L-S\equiv (d-1)S$. 
Thus, $\wt{L}-F=L'-S'=\pi^{*}(L-S)$ is nef. 

From now on, we put the index $k$ to divisors and varieties, 
i.e. $V,D,S,\wt{D}$ etc. are written as $V_{k},D_{k},S_{k},\wt{D}_{k}$ etc. 
Recall that $V_{k}$ is a smooth hypersurface of degree $d\geq 1$ 
in $\PP^{k}$ and $D_{k}=\PP^{m}\times V_{k}$. 
Assume $[\wt{L}_{k}-F_{k}]\in \Nef (\wt{D}_{k})$. 
Consider a smooth divisor $M_{k+1}$ on $D_{k+1}=\PP^{m}\times V_{k+1}$ 
such that $M_{k+1}\sim L_{k+1}$ and $S_{k+1}\subset M_{k+1}$. 
Put $\wt{M}_{k+1}:=(\pi\circ \beta)\st M_{k+1}$. 
Since $M_{k+1}$ is isomorphic to $\PP^{m}\times V_{k}=D_{k}$ 
as variety, so is 
$\wt{M}_{k+1}$ to $\wt{D}_{k}$. 
The identification $\wt{M}_{k+1}=\wt{D}_{k}$ yields
\begin{eqnarray*}
\wt{L}_{k+1}|_{\wt{M}_{k+1}}=\wt{L}_{k},\ \ \  
F_{k+1}|_{\wt{M}_{k+1}}=F_{k}. 
\end{eqnarray*} 
Hence, by inductive hypothesis, 
$[\wt{L}_{k+1}-F_{k+1}]|_{\wt{M}_{k+1}}\in \Nef (\wt{M}_{k+1})$.
Since $\wt{M}_{k+1}\sim \wt{L}_{k+1}-F_{k+1}$, we have 
\begin{eqnarray*}
[\wt{L}_{k+1}-F_{k+1}-\wt{M}_{k+1}]=0 
\in \Nef(\wt{D}_{k+1}).
\end{eqnarray*}
By Lemma \ref{induction}, $[\wt{L}_{k+1}-F_{k+1}]\in \Nef (\wt{D}_{k+1})$. 
Therefore, $[\wt{L_{k}}-F_{k}]\in \Nef(\wt{D}_{k})$ for any $k\geq 2$.  
\qed

\section{Structure of nef cone and effective cone}\label{cone}
From now on, we fix the following.

Notation (*): 
Let $n$ and $k$ be natural numbers such that $n\geq 3$ and  $2\leq k\leq n-1$. 

Put $Y:=\PP^{n-k}\times \PP^{k}$.  
Let $p:Y\to \PP^{n-k}$ and $q:Y\to \PP^{k}$ be the two projections. 
Put $H:=p^{*}\OO_{\PP^{n-k}}(1)$ and $L:=q^{*}\OO_{\PP^{k}}(1)$. 
Let $P_{0}$ be a fiber of $p:Y\to \PP^{n-k}$. 
Let $C$ be a smooth hypersurface of degree $d\geq 1$ in $P_{0}\simeq \PP^{k}$.
Let $S$ be a fiber of $q:Y\to \PP^{k}$ such that $S\cap C\neq \emptyset$. 
Put $y_{0}:=S\cap C =S\cap P_{0}$. 
Let $h$ be a line in a fiber of $p$ such that $h\cap S=\emptyset$ and $h\cap C=\emptyset$. 
Let $l$ be a line in a fiber of $q$ such that $l\cap S=\emptyset$ and $l\cap C=\emptyset$. 
Let $h_{0}$ be a line in $P_{0}\simeq \PP^{k}$ such that $h_{0}\not\subset C$ and $y_{0}\not\in h_{0}$. 
Let $l_{0}$ be a line in a fiber of $q$ different from $S$ such that $l_{0}\cap C\neq \emptyset$. 
Let $H_{0}$ be an element of the the linear system $|H|$ such that $P_{0}\subset H_{0}$. 
Let $L_{0}$ be an element of the linear system $|L|$ such that $S\subset L_{0}$. If $d=1$, we assume that $L_{0}$ does not contain $C$. 
Put $D:=q^{-1}(q(C))$. Note that we have $D\sim dL$. 

Let $\pi:X\to Y$ be the blow-up along $C$. 
Put $E:=\Exc (\pi)$ and $E_{0}:=\pi^{-1}(y_{0})$. 
Let $e_{0}$ be a line in $E_{0}\simeq \PP^{n-k}$. 
Let $e$ be a line in a fiber different from $E_{0}$ 
of the $\PP^{n-k}$-bundle $\pi |_{E}:E\to C$ 
(if $k=n-1$, then $e$ and $e_{0}$ 
are fibers of the $\PP^{1}$-bundle $\pi |_{E}$). 
Put $S':=\pi\st S$ and $P_{0}':=\pi\st P_{0}$. 
Note that 
we have $S\cap P_{0}=y_{0}$ while $S'\cap P_{0}'=\emptyset$.  
Put $H':=\pi^{*}H$ and $L':=\pi^{*}L$. 
Let $H_{0}',L_{0}',D',h',l',h_{0}'$ and $l_{0}'$ 
be the strict transforms by $\pi$ of 
$H_{0},L_{0},D,h,l,h_{0}$ and $l_{0}$, respectively.

Let $\beta:\wt{X}\to X$ be the blow-up along $S'$. 
Put $F:=\Exc (\beta)$. 
Let $f$ be a line in a fiber of the $\PP^{k-1}$-bundle $\beta |_{F}:F\to S'$  
(if $k=2$ then $f$ is a fiber of the $\PP^{1}$-bundle $\beta |_{F}$). 
Put 
$\wt{H}:=\beta^{*}H'$,
$\wt{L}:=\beta^{*}L'$ 
and $\wt{E}:=\pi^{*}E$. 
Finally, let 
$\wt{H_{0}},\wt{L_{0}},\wt{D},
\wt{h},\wt{l},\wt{h_{0}},\wt{l_{0}},\wt{e}$ and $\wt{e_{0}}$ 
be the strict transforms by $\beta$ of 
$H_{0}',L_{0}',D',h',l',h_{0}',l_{0}',e$ and $e_{0}$, respectively. 

\begin{lemma}\label{intersection}
We have the following table of intersection numbers. 
\begin{eqnarray*}
\begin{array}{c|cccc}
&\wt{H}&\wt{L}&\wt{H}+d\wt{L}-\wt{E}&\wt{H}+d\wt{L}-\wt{E}-F \\ 
\hline
\wt{l_{0}}&1&0&0&0 \\ 
\wt{h_{0}}&0&1&0&0 \\ 
\wt{e_{0}}&0&0&1&0 \\ 
f&0&0&0&1
\end{array}
\end{eqnarray*}
\end{lemma}
{\it Proof}. 
We have $e\equiv e_{0}$ in $X$. Note that $e_{0}$ intersects $S'$ transversally at one point. 
Hence $\wt{e}\equiv \wt{e}_{0}+f$ in $\wt{X}$. 
We have $l \equiv l_{0}$ in $Y$. 
Since $l_{0}$ intersects $C$ transversally at one point, 
we have $l'\equiv l_{0}'+e$ in $X$. 
This implies $\wt{l}\equiv \wt{l_{0}}+\wt{e}$ in $\wt{X}$. 
Note that $h\equiv h_{0}$ and $(C\cdot h_{0})_{P_{0}}=d$. 
Hence $h'\equiv h'_{0}+de$, i.e. $\wt{h}\equiv \wt{h_{0}}+d\wt{e}$ 
 because $S' \cap h_{0}'=\emptyset$.  

The following intersection numbers are obvious.
\begin{eqnarray*}
\begin{array}{c|cccc}
&\wt{H}&\wt{L}&\wt{E}&F \\ 
\hline
\wt{l} &1&0&0&0 \\
\wt{h} &0&1&0&0 \\ 
\wt{e} &0&0&-1&0 \\ 
f &0&0&0&-1
\end{array}
\end{eqnarray*}
Therefore, we have the following table. 
\begin{eqnarray*}
\begin{array}{c|cccc|cccc}
&\wt{H}&\wt{L}&\wt{E}&F &\wt{H}&\wt{L}&\wt{H}+d\wt{L}-\wt{E}&\wt{H}+d\wt{L}-\wt{E}-F \\ 
\hline
\wt{l_{0}}\equiv \wt{l}-\wt{e} &1&0&1&0&1&0&0&0 \\ 
\wt{h_{0}}\equiv \wt{h}-d\wt{e} &0&1&d&0&0&1&0&0 \\ 
\wt{e_{0}}\equiv \wt{e}-f &0&0&-1&1&0&0&1&0 \\ 
f &0&0&0&-1&0&0&0&1
\end{array}
\end{eqnarray*}
\qed

\begin{proposition}\label{nefcone} We have 
\begin{eqnarray*}
\Nef (\wt{X})
=\RR [
\wt{H},\wt{L},\wt{H}+d\wt{L}-\wt{E}, 
\wt{H}+d\wt{L}-\wt{E}-F].
\end{eqnarray*}
\end{proposition}
{\it Proof}. We have $[\wt{H}],[\wt{L}]\in\Nef(\wt{X})$ and 
$[\wt{l_{0}}],[\wt{h_{0}}],[\wt{e_{0}}],[f]\in \NNE(\wt{X})$. 
Recall $\rho (\wt{X})=4$. 
By Lemma \ref{intersection}, it is sufficient to show that 
$\wt{H}+d\wt{L}-\wt{E}$ and $\wt{H}+d\wt{L}-\wt{E}-F$ are nef (see \cite{Tsukioka} Lemma 3.1). 

First, we show $[\wt{H}+d\wt{L}-\wt{E}-F]\in\Nef (\wt{X})$. 
Put $V:=q(C)$. We have $D\simeq \PP^{n-k}\times V$. Consider the projections
$pr_{1}:=p|_{D}:D\to \PP^{n-k}$ and $pr_{2}:=q|_{D}:D\to V$. Note that 
$C$ is a fiber of $pr_{1}$ and $S$ is a fiber of $pr_{2}$. 
If we put $m:=n-k$, 
then our divisors $\wt{H}|_{\wt{D}},\wt{L}|_{\wt{D}},\wt{E}|_{\wt{D}},F|_{\wt{D}}$ 
correspond exactly to the divisors $\wt{H},\wt{L},\wt{E},F$ in Proposition \ref{nef}. 
Hence $[\wt{H}+d\wt{L}-\wt{E}-F]|_{\wt{D}}\in\Nef (\wt{D})$. 
Note that $\wt{D}\sim d\wt{L}-\wt{E}-F$. 
By Lemma \ref{induction}, 
we conclude  
$[\wt{H}+d\wt{L}-\wt{E}-F]\in \Nef (\wt{X})$. 

Now, we show $[\wt{H}+d\wt{L}-\wt{E}]\in\Nef(\wt{X})$. 
If $n-k=1$, we have $S'\simeq S\simeq \PP^{1}$, 
$H'|_{S'}\simeq \OO_{\PP^{1}}(1)$ and $E|_{S'}\simeq \OO_{\PP^{1}}(1)$. 
Hence $H'|_{S'}-E|_{S'}$ is linearly equivalent to $0$, in particular nef. 
If $n-k\geq 2$, then $\pi |_{S'}:S'\to S=\PP^{n-k}$ is the blow-up at the point $y_{0}$. 
By Lemma \ref{point}, $H'|_{S'}-E|_{S'}=(\pi |_{S'})^{*}\OO_{\PP^{n-k}}(1)-\Exc (\pi |_{S'})$ is nef. 
Since $L|_{S}\sim 0$, we have $L'|_{S'}\sim 0$. 
Hence, 
$[H'+dL'-E]|_{S'}=[H'-E]|_{S'}\in\Nef(S')$, 
which implies 
$$[\wt{H}+d\wt{L}-\wt{E}]|_{F}
=(\beta |_{F})^{*}\left([H'+dL'-E]|_{S'}\right)\in\Nef (F).$$
Since $[\wt{H}+d\wt{L}-\wt{E}-F]\in\Nef(\wt{X})$ 
as shown above, 
we have by Lemma \ref{induction} 
$[\wt{H}+d\wt{L}-\wt{E}]\in \Nef(\wt{X})$. 
 \qed

\

Remark: $H'-E$ itself is not nef because $(H'-E)|_{P_{0}'}\simeq \OO_{\PP^{k}}(-d)$.

\begin{proposition}\label{effective cone} We have 
\begin{eqnarray*}
\Eff (\wt{X})
=\RR [\wt{H_{0}},\wt{L_{0}},\wt{E},F, \wt{D}]. 
\end{eqnarray*}
\end{proposition}
{\it Proof}. The divisors $\wt{H_{0}},\wt{L_{0}},\wt{E},F,\wt{D}$ are all effective. 
Hence, the inclusion 
$\RR [\wt{H_{0}},\wt{L_{0}},\wt{E},F, \wt{D}] \subset \Eff(\wt{X})$ holds. 
Let $A\subset \wt{X}$ be a prime divisor. Suppose $A\not\in\{ \wt{E},F, \wt{D}\}$. 
Since $B:=(\pi\circ \beta)(A)$ is a prime divisor on $Y=\PP^{n-k}\times \PP^{k}$, 
there exist $a,b\in \Z_{\geq 0}$ such that $B\sim aH+bL$. 
Put $\mu:=\mult_{C}B$ and $\nu:=\mult_{S}B$. Since $A=(\pi\circ\beta)\st B$, we have 
$A\sim a\wt{H}+b\wt{L}-\mu\wt{E}-\nu F$. 
If $q^{-1}(u)\subset B$ for any $u\in q(C)$, 
we have $D=q^{-1}(q(C)))\subset B$, which implies $B=D$, a contradiction to the choice of $A$. 
Hence, there exists a point $u\in q(C)$ 
such that $q^{-1}(u)\not\subset B$. 
Put $Y_{u}:=q^{-1}(u)$ and $c_{u}:=C\cap Y_{u}$. 
Note that there exists a line 
$l_{u} \subset Y_{u}\simeq \PP^{n-k}$ 
such that $c_{u}\in l_{u}$ and $l_{u}\not\subset B$ 
(if $n-k=1$, we put $l_{u}:=Y_{u}\simeq \PP^{1}$). 
Since $B\sim aH+bL$ and 
$l_{u}\equiv l$, we have $B\cdot l_{u}=a$  
(if $n-k=1$,  
we have $B |_{Y_{u}}\cdot l_{u}=\deg (B |_{Y_{u}})$). 
Therefore, 
$$
a=B\cdot l_{u}=B |_{Y_{u}}\cdot l_{u}\geq \mult_{c_{u}}(B |_{Y_{u}})\geq \mult_{C}B=\nu.
$$
Since $B$ is a prime divisor on $Y=\PP^{n-k}\times \PP^{k}$, there exists a point $v\in \PP^{n-k}$ 
such that $p^{-1}(v)\not\subset B$. Put $P_{v}:=p^{-1}(v)$ and $s_{v}:=S\cap P_{v}$. 
Note that there exists a line $h_{v}\subset P_{v}\simeq \PP^{k}$ such that 
$s_{v}\in h_{v}$ and $h_{v}\not\subset B$. Since $B\sim aH+bL$ and $h_{v}\equiv h$, we have 
$B\cdot h_{v}=b$. Hence, 
$$
b=B\cdot h_{v}=B |_{P_{v}}\cdot h_{v}\geq \mult_{s_{v}}(B |_{P_{v}})\geq \mult_{S}B=\nu.
$$ 
Recall that $\wt{H}\sim \wt{H_{0}}+\wt{E}$ and $\wt{L}\sim \wt{L_{0}}+F$. 
We conclude 
$$
[A]=[a\wt{H}+b\wt{L}-\mu\wt{E}-\nu F] 
=[a\wt{H_{0}}+b\wt{L_{0}}+(a-\mu)\wt{E}+(b-\nu)F]
\in \RR[\wt{H_{0}},\wt{L_{0}},\wt{E},F].
$$
\qed 

Remark: 
Since $\wt{D}\sim d\wt{L}-\wt{E}-F\sim d\wt{L_{0}}-\wt{E}+(d-1)F$, 
we have 
$[\wt{D}]\not\in\RR[\wt{H_{0}},\wt{L_{0}},\wt{E},F]$ 
for any $d\geq 2$. 
If $d=1$, 
we have 
$\wt{L_{0}}\sim \wt{D}+\wt{E}$ and 
$\RR[\wt{H_{0}},\wt{L_{0}},\wt{E},F,\wt{D}]=\RR[\wt{H_{0}},\wt{E},F,\wt{D}]$. 

\section{Proof of Theorem}\label{proof}

We continue to use the notation (*). First, we find conditions for $-(K_{\wt{X}}+\Delta)$ to be ample 
where $\Delta$ is a $\Q$-divisor on $\wt{X}$. 

\begin{lemma}\label{inequalities} 
Consider a divisor $\Delta=\alpha\wt{H}+\beta\wt{L}+\gamma\wt{E}+\delta F$ with $\alpha,\beta,\gamma,\delta\in\Q$. 
Then, $-(K_{\wt{X}}+\Delta)$ is ample if and only if 
\begin{eqnarray*}
\begin{cases}
\alpha+\gamma  &<1 \\ 
\beta+d\gamma &<k+1-(n-k)d \\ 
-\gamma+\delta &<n-2k+1 \\ 
-\delta &<k-1.
\end{cases}
\end{eqnarray*}
\end{lemma}
{\it Proof}. Since $Y=\PP^{n-k}\times \PP^{k}$, we have 
$-K_{Y}\sim (n-k+1)H+(k+1)L$. 
By the adjunction formula for the blow-ups $\pi$ and $\beta$, 
\begin{eqnarray*}
-K_{\wt{X}}
&\sim& \beta^{*}(-K_{X})-(k-1)F \\ 
&\sim& \beta^{*}(\pi^{*}(-K_{Y})-(n-k)E)-(k-1)F \\ 
&\sim& (n-k+1)\wt{H}+(k+1)\wt{L}-(n-k)\wt{E}-(k-1)F.
\end{eqnarray*}
Hence, $-(K_{\wt{X}}+\Delta)$ is linearly equivalent to 
\begin{eqnarray*}
&&(n-k+1-\alpha)\wt{H}+(k+1-\beta)\wt{L}-(n-k+\gamma)\wt{E}-(k-1+\delta)F \\ 
&=&(1-\alpha-\gamma)\wt{H}+(k+1+dk-dn-\beta-d\gamma)\wt{L} \\ 
&&+(n-2k+1+\gamma-\delta)(\wt{H}+d\wt{L}-\wt{E})
+(k-1+\delta)(\wt{H}+d\wt{L}-\wt{E}-F).
\end{eqnarray*} 
By Proposition \ref{nefcone},  we conclude that $-(K_{\wt{X}}+\Delta)$ is ample (i.e. the numerical equivalence class is an interior point of the nef cone) if and only if 
\begin{eqnarray*}
\begin{cases}
1-\alpha-\gamma &> 0 \\ 
k+1+dk-dn-\beta-d\gamma &>0 \\ 
n-2k+1+\gamma -\delta &>0 \\ 
k-1+\delta &>0.
\end{cases}
\end{eqnarray*}
\qed

Remark: We obtain the condition for $-(K_{\wt{X}}+\Delta)$ to be nef if we replace ``$<$" by ``$\leq$". 

\begin{lemma}\label{ample nef} 
$-K_{\wt{X}}$ is ample (resp. nef) if and only if 
\begin{eqnarray*}
\begin{cases}
0 &< k+1-(n-k)d \\ 
0 &< n-2k+1 
\end{cases} \ \ 
\left( \mbox{resp.} 
\begin{cases}
0 &\leq k+1-(n-k)d \\ 
0 &\leq n-2k+1 
\end{cases}
\right).
\end{eqnarray*}
\end{lemma}
{\it Proof}. We put $\Delta=0$ in Lemma \ref{inequalities}. \qed

\

In view of Lemma \ref{effective cone}, it is natural that we should consider boundary divisors 
in the form: 
$$\Delta=x\wt{H_{0}}+y\wt{L_{0}}+z\wt{E}+wF+u\wt{D} \ \ \ (x,y,z,w,u\in\Q^{+}).$$

\begin{lemma}\label{log Fano pair} 
$(\wt{X},\Delta)$ is a log Fano pair if and only if 
\begin{eqnarray*}
\begin{cases}
0\leq x,y,z,w,u &<1\\ 
z-u &<1 \\ 
y-dx+dz &<k+1-(n-k)d \\ 
x-y-z+w &<n-2k+1 \\ 
y-w+u &<k-1.
\end{cases}
\end{eqnarray*}
\end{lemma}
{\it Proof}. We see that $\Supp\Delta=\wt{H_{0}}\cup\wt{L_{0}}\cup\wt{E}\cup F\cup\wt{D}$ 
is a simple normal crossing divisor. Hence, the pair $(\wt{X},\Delta)$ is klt if and only if 
$0\leq x,y,z,w,u<1$. 
On the other hand,  
\begin{eqnarray*}
\Delta
&\sim&x(\wt{H}-\wt{E})+y(\wt{L}-F)+z\wt{E}+wF+u(d\wt{L}-\wt{E}-F) \\ 
&=&x\wt{H}+(y+du)\wt{L}+(-x+z-u)\wt{E}+(-y+w-u)F. 
\end{eqnarray*} 
If we put $\alpha:=x$, $\beta:=y+du$, $\gamma:=-x+z-u$ and $\delta:=-y+w-u$, 
then, 
$\alpha+\gamma=z-u$, $\beta+d\gamma=y-dx+dz$, 
$-\gamma+\delta=x-y-z+w$ and $-\delta=y-w+u$. Thus, the statement follows from Lemma \ref{inequalities}. \qed

\

Now, we prove Theorem, which is divided into the following Propositions. 

\begin{proposition}\label{big} The anti-canonical divisor $-K_{\wt{X}}$ is big for any $n\geq 3$, $k\in \{ 2,\cdots, n-1\}$ and $d\geq 1$. 
\end{proposition}
{\it Proof}. Recall $-K_{\wt{X}}\sim (n-k+1)\wt{H}+(k+1)\wt{L}-(n-k)\wt{E}-(k-1)F$. 
We consider the divisors 
\begin{eqnarray*}
A&:=& \left(\frac{1}{2}+\frac{1}{d}\right)\wt{H}
+\left(1+\frac{1}{2d}\right)\wt{L}-\frac{1}{d}\wt{E}-\frac{1}{2d}F, \\ 
B&:=& \left(n-k+\frac{1}{2}-\frac{1}{d}\right)\wt{H}
+\left(k-\frac{1}{2d}\right)\wt{L}
-\left(n-k-\frac{1}{d}\right)\wt{E}-\left(k-1-\frac{1}{2d}\right)F.
\end{eqnarray*}
Then, we have $-K_{\wt{X}}\sim A+B$. 
On the other hand, 
\begin{eqnarray*}
A&\sim& \frac{1}{2}\wt{H}+\frac{1}{2d}\wt{L}
+\frac{1}{2d}\left(\wt{H}+d\wt{L}-\wt{E}\right)
+\frac{1}{2d}\left(\wt{H}+d\wt{L}-\wt{E}-F\right), \\
B&\sim& \left(n-k+\frac{1}{2}-\frac{1}{d}\right)\wt{H_{0}}
+\left(k-\frac{1}{2d}\right)\wt{L_{0}}+\frac{1}{2}\wt{E}+F.
\end{eqnarray*} 
By Proposition \ref{nefcone}, $A$ is ample. 
Since $n-1\geq k\geq 2$, we have $n-k+\frac{1}{2}-\frac{1}{d}>0$ 
and $k-\frac{1}{2d}>0$. 
Hence, $B$ is effective. The statement follows from  \cite{Lazarsfeldbook} Corollary 2.2.7. \qed

\

Remark: It follows that $\wt{X}$ is weak Fano if and only if 
$-K_{\wt{X}}$ is nef. 

\begin{proposition}\label{n=3}
Assume $n=3$ (hence $k=2$). Then, $\wt{X}$ is log Fano for any $d\geq 1$. Moreover 
$\wt{X}$ is weak Fano if and only if $d\in \{1,2,3\}$, and is not Fano for any $d\geq 1$. 
\end{proposition}
{\it Proof}.  If $d=1$, we put $z=\frac{1}{2}$ and $x=y=w=u=0$. 
If $d\geq 2$, we put 
\begin{eqnarray*}
x=y=\frac{d-2}{d-1}, \ 
z=\frac{1}{2d} \ 
\mbox{ and } \ 
w=u=0. 
\end{eqnarray*} 
Then, the inequalities in Lemma \ref{log Fano pair} are all satisfied. 
Hence, $\wt{X}$ is log Fano for any $d\geq 1$. 
By Lemma \ref{ample nef}, 
$-K_{\wt{X}}$ is nef if and only if $d\in \{1,2,3\}$, and is not ample 
for any $d\geq 1$. \qed

\

From now on, we assume $n\geq 4$. 

\begin{proposition}\label{weakFano}
$\wt{X}$ is weak Fano if and only if 
\begin{eqnarray*}
(n,k,d)=(2k-1,k,1),(2k,k,1),(2k+1,k,1) \ \mbox{or} \ (5,3,2).
\end{eqnarray*}
$\wt{X}$ is Fano if and only if  $(n,k,d)=(2k,k,1)$.
\end{proposition}
{\it Proof}. By Lemma \ref{ample nef} and the remark to Proposition \ref{big},  
$\wt{X}$ is weak Fano if and only if 
\begin{eqnarray*}
\begin{cases}
dn & \leq  (d+1)k+1 \\ 
2k-1 & \leq  n. 
\end{cases}
\end{eqnarray*}
This yields $d(2k-1)\leq dn \leq (d+1)k+1.$
In particular, we have $d(2k-1)\leq (d+1)k+1$, which implies 
\begin{eqnarray*}
d\leq \frac{k+1}{k-1}=1+\frac{2}{k-1}\leq 1+2=3. 
\end{eqnarray*}

If $d=1$, we have $2k-1\leq n \leq 2k+1$.
If $d=2$, we have $4k-2\leq 2n \leq 3k+1$. 
Since $n\geq 4$, we have $8\leq 3k+1$, i.e. $k\geq 3$. 
On the other hand, we have $4k-2\leq 3k+1$, i.e. $k\leq 3$. 
Thus, $k=3$ and $n=5$. 
If $d=3$, we have $6k-3\leq 3n\leq 4k+1$. 
Since $n\geq 4$, we have $12\leq 4k+1$, i.e. $k\geq 3$. 
On the other hand, we have $6k-3\leq 4k+1$, i.e. $k\leq 2$, a contradiction. 
Therefore, we have 
\begin{eqnarray*}
(n,k,d)=(2k-1,k,1),(2k,k,1),(2k+1,k,1) \ \mbox{or} \ (5,3,2).
\end{eqnarray*}

In any of these cases, $k+1-(n-k)d$ and $n-2k+1$ are non-negative, and strictly positive only for the case $(n,k,d)=(2k,k,1)$. \qed 

\begin{proposition}\label{log Fano} 
There exist rational numbers $x,y,z,w,u$ such that 
$(\wt{X},\Delta=x\wt{H_{0}}+y\wt{L_{0}}+z\wt{E}+wF+u\wt{D})$ is a log Fano pair 
if and only if one of the following holds: 
\begin{itemize}
\item $n\in\{2k-2,2k-1,2k,2k+1\}$ and $d=1$, 
\item $(n,k,d)=(4,2,2),(4,3,2) \mbox{ or } (5,3,2)$.
\end{itemize}
\end{proposition}
{\it Proof}. Assume that $(\wt{X},\Delta)$ is a log Fano pair. 
By Lemma \ref{log Fano pair}, we have $0\leq x,y,z,w<1$ and 
\begin{eqnarray}\label{logfano1}
\begin{cases}
y-dx+dz &<k+1-(n-k)d \\ 
x-y-z+w&<n-2k+1.
\end{cases}
\end{eqnarray}
We have 
$-d<-dx\leq y-dx+dz<k+1-nd+kd$, i.e. $nd<(k+1)(d+1)$. 
On the other hand, $-2<-y-z\leq x-y-z+w<n-2k+1$, i.e. $2k-3<n$. 
Thus, 
\begin{eqnarray}\label{logfano2}
2k-3<n<\frac{(k+1)(d+1)}{d}.
\end{eqnarray}
If $d=1$, we have $2k-3<n<2k+2$, i.e. $n\in\{ 2k-2,2k-1,2k,2k+1\}$.

Consider the case $d\geq 2$. Since $w\geq 0$, 
we have $x-y-z\leq x-y-z+w$.  
Hence, by (\ref{logfano1}),  we obtain:  
\begin{eqnarray*}
\begin{cases}
y-d(x-z) &<k+1-nd+kd \\ 
-y+(x-z) &<n-2k+1.
\end{cases}
\end{eqnarray*}
The sum of two inequalities gives $(d-1)(n-k)-2<(d-1)(x-z)$, 
i.e. $n-k-\frac{2}{d-1}<x-z$. 
On the other hand, we have 
$x-z<n-2k+1+y<n-2k+2$. Therefore, $n-k+\frac{2}{1-d}<n-2k+2$, i.e. 
\begin{eqnarray}\label{logfano3}
k<2+\frac{2}{d-1}. 
\end{eqnarray}
Since $d,k\geq 2$, this implies $k=2$ or $3$. 
If $k=2$, by (\ref{logfano2}), we have $n<\frac{3(d+1)}{d}$. Thus, $n=4$ and $d=2$. 
If $k=3$, by (\ref{logfano3}), we have $d=2$. This implies $3<n<6$ 
by (\ref{logfano2}). Thus, $n=4$ or $5$.
Therefore, we have 
\begin{eqnarray*}
(n,k,d)=(4,2,2),(4,3,2) \mbox{ or } (5,3,2). 
\end{eqnarray*}

In each case, we put $(x,y,z,w,u)$ as follows:  
\begin{center}
\begin{tabular}{l|c}
\ \ \ \ \ \ $(n,k,d)$ & $(x,y,z,w,u)$ \\
\hline 
Case: $n\geq 4, d=1$ & \\ 
\ \ $n=2k-2$ & $(0,\frac{1}{2},\frac{3}{4},0,0)$ \\ 
\ \ $n=2k-1$ & $(0,\frac{1}{2},0,0,0)$ \\ 
\ \ $n=2k$ & $(0,0,0,0,0)$ \\ 
\ \ $n=2k+1$ & $(\frac{1}{2},0,0,0,0)$ \\  
\hline 
$(n,k,d)=(4,2,2)$ & $(\frac{3}{4},0,0,0,0)$\\ 
$(n,k,d)=(4,3,2)$ & $(\frac{1}{8},\frac{1}{2},\frac{3}{4},0,0)$\\ 
$(n,k,d)=(5,3,2)$ & $(\frac{1}{2},\frac{1}{2},\frac{1}{8},0,0)$
\end{tabular}
\end{center}
Then, the inequalities in Lemma \ref{log Fano pair} are satisfied. \qed

------------

{\small {\sc Department of Mathematics, 
Tokai University, 
4-1-1 Kitakaname, Hiratsuka-shi,
Kanagawa, 259-1292 Japan} }

{\it E-mail address}: \ {\tt tsukioka@tokai-u.jp}
 
\end{document}